\documentclass[12pt]{article}
\usepackage{amsfonts,amsmath,amsthm,amssymb,xypic,multicol}

\DeclareMathOperator{\Gal}{{Gal}}
\DeclareMathOperator{\md}{{mod}}

\DeclareMathOperator{\chr}{{char}}

\newtheorem{defn}{Definition}[section]
\newtheorem{prop}[defn]{Proposition}
\newtheorem{thm}[defn]{Theorem}
\newtheorem{cor}[defn]{Corollary}

\newtheorem{lemma}[defn]{Lemma}

\newcommand{\R}{\ensuremath{\mathbb{R}}}
\newcommand{\p}{\ensuremath{\mathfrak{p}}}
\newcommand{\q}{\ensuremath{\mathfrak{q}}}

\newcommand{\F}{\ensuremath{\mathbb{F}}}
\newcommand{\Z}{\ensuremath{\mathbb{Z}}}

\newcommand{\Q}{\ensuremath{\mathbb{Q}}}
\newcommand{\OK}{\ensuremath{\mathcal{O}}}
\newcommand{\w}{\omega}

\title{Function Fields with Class Number Indivisible by a Prime $\ell$}
\author{Michael Daub, Jaclyn Lang, Mona Merling, Allison M. Pacelli,\\  Natee Pitiwan, and Michael Rosen}

\begin{document}
\maketitle


\begin{abstract}
It is known that infinitely many number fields and function fields of any degree $m$ have class number divisible by a given integer $n$.  However, significantly less is known about the indivisibility of class numbers of such fields.  While it's known that there exist infinitely many quadratic number fields with class number indivisible by a given prime,
the fields are not constructed explicitly, and nothing appears to be known for higher degree extensions.  In \cite{Pacelli-Rosen}, Pacelli and Rosen explicitly constructed an infinite class of function fields of any degree $m$, $3 \nmid m$, over $\F_q(T)$ with class number indivisible by $3$, generalizing a result of Ichimura for quadratic extensions.  Here we generalize that result, constructing, for an arbitrary prime $\ell$, and positive integer $m > 1$, infinitely many function fields of degree $m$ over the rational function field, with class number indivisible by $\ell$.
\end{abstract}

\section{Introduction}
The question of class number indivisibility has always been more difficult than the question of class number divisibility.  For example, although Kummer was able to prove Fermat's Last
Theorem for regular primes, that is, primes $p$ not dividing the class
number of the $p$-th cyclotomic field, it is still unknown today whether
infinitely many regular primes exist (in 1915, Jensen did prove the
existence of infinitely many irregular primes).

In 1976, Hartung~\cite{Ha} showed that infinitely many imaginary quadratic
number fields have class number not divisible by 3.  The analogous result
for function fields was proven in 1999 by Ichimura~\cite{Ichimura}.  Horie and
Onishi \cite{Ho, Ho2, HoO}, Jochnowitz~\cite{Jo}, and Ono and
Skinner~\cite{OS} proved that there are infinitely many imaginary
quadratic number fields with class number not divisible by a given prime
$p$. Quantitative results on the density of quadratic fields with class
number indivisible by 3 have been obtained by Davenport and
Heilbronn~\cite{DH}, Datskovsky and Wright~\cite{DW}, and Kimura~\cite{Ki}
(for relative class numbers). Kohnen and Ono made further progress
in~\cite{KO}.  They proved that for all $\epsilon > 0$ and sufficiently
large $x$,  the number of imaginary quadratic number fields $K =
\Q(\sqrt{-D})$ with $p \nmid h_K$ and $D < x$ is $$\geq
\left(\frac{2(p-2)}{\sqrt3(p-1)} -
\epsilon\right)\frac{\sqrt{x}}{\log{x}}.$$   Less is known about
class numbers in real quadratic fields, but in 1999, Ono~\cite{O} obtained
a similar lower bound for the number of real quadratic fields $K$ with $p
\nmid h_K$ and bounded discriminant; this bound is valid for primes $p$
with  $3 < p < 5000$.  The results above do not give explicit families of
fields with the desired class number properties.  In 1999,
Ichimura~\cite{Ichimura} constructed an explicit infinite family of quadratic
function fields with class number not divisible by 3.  Pacelli and Rosen~\cite{Pacelli-Rosen} extended this to non-quadratic fields of degree $m$ over $\F_q(T)$, $3 \nmid m$.
In this paper, we
generalize Pacelli and Rosen's result, constructing, for a large class of $q$, infinitely many function fields of any degree $m$ over $\F_q(T)$ with class number indivisible by an arbitrary prime $\ell$.

For similar results on divisibility of class numbers,  see
Nagell~\cite{Nag} for imaginary number fields, Yamamoto~\cite{Y} or
Weinberger~\cite{W} for real number fields, and Friesen~\cite{F} for
function fields.  For quantitative results, see for example Murty~\cite{murty1, murty2}.  More generally, to see results on the minimum $n$-rank of the ideal class group of a global field, see Azuhata and Ichimura~\cite{AI} or Nakano~\cite{Nak}
for number fields and Lee and Pacelli~\cite{Lee2, LP1, LP2, P1, P2} for function fields.

As in~\cite{Pacelli-Rosen}, the fields we construct are given explicitly.  The idea of the
proof is to construct two towers of fields $N_1 \subset \cdots \subset N_t
= \F_q(T)$ and ${M_1 \subset \cdots \subset M_t} $.  The fields are designed
so that  $\ell \nmid h_{M_1}$, $N_{i+1}/N_i$ is cyclic of degree $\ell$ and ramified (totally)
at exactly one prime, $M_i/N_i$ is a degree $m$ extension, and $M_{i+1}$
is the composite field of $M_i$ and $N_{i+1}$.  Together with class field
theory, this is enough to show that $\ell \nmid h_{M_i}$ for any $1 < i \leq
t$.  Thus $M_t$ has  degree $m$ over $N_t$, the rational function field,
and has class number not divisible by $\ell$.

Let $q$ be a power of an odd prime, and $\F_q$ the finite field with $q$ elements. 
The main results are as follows:

\begin{thm}\label{t:main_result}
Let $m$ be any positive integer $m > 1$ and $\ell$ an odd prime. Write $m = \ell^tm_1$ for integers $t$ and $m_1$ with $\ell \nmid m_1$.  Let $m_0$ be the square-free part of $m_1$, and assume that $q$ is sufficiently large with $q \equiv 1$ (mod $m_0$) and $q \equiv -1$ (mod $\ell$).  Then there are infinitely many function fields $K$ of degree $m$ over  $\F_q(T)$ with $\ell \nmid h_K$.
\end{thm}

\begin{cor}\label{c_cor1}
If $q$ satisfies the hypotheses of the theorem and, in addition, if $q \equiv 1$ (mod $m$), then there are infinitely many cyclic extensions $K$ of degree $m$ over $\F_q(T)$ with $\ell \nmid h_K$.
\end{cor}

\begin{cor}\label{c:cor2}
If $q$ satisfies the hypotheses of the theorem and, in addition, $m$ is square-free and $q \equiv 1$ (mod $m_1$), then there are infinitely many cyclic extensions $K$ of degree $m$ over $\F_q(T)$ with $\ell \nmid h_K$.
\end{cor}

For the remainder of this introduction, we will outline some important results and methods which will be used in the proof of the main theorem, Theorem 1.1, above.  In the statement of Theorem 1.1 we use the phrase ``all sufficiently large $q$.''  In the Appendix we will give a quantitative version of this restriction.  In Section 3, we prove a function field analogue of a class field theoretic result of Iwasawa; this result is stated but not proved by Ichimura in \cite{Ichimura}.  In Section 4, we prove the main theorem, and in Section 5, we prove the two corollaries stated above.

In [27] the cubic extensions needed were generated by using a variant of the ``simplest cubic polynomials" discovered by Dan Shanks [30]; $ X^3-3uX^2-(3u+3)X-1$. Any root of this polynomial generates a Galois extension of $k(u)$ with Galois group isomorphic to $\mathbb{Z}/3 \mathbb{Z} $. Here $k$ is any field with characteristic different from $3$. Hashimoto and Miyake found generalizations of this polynomial for any odd degree $\ell$. Their work was simplified and extended by Rikuna in [28] and further developed by Komatsu in [17]. We will restrict ourselves to the case $\ell$ is odd and present Rikuna's polynomials as exposited in Komatsu.\par

Let $K$ be a field whose characteristic does not divide $\ell$. Let $\zeta$ be a primitve $\ell$-th root of unity in some field containing $k$ and suppose $\omega=\zeta+\zeta^{-1}$ is in $K$. Define

$$              {\cal P}(X) := (\zeta^{-1}-\zeta)^{-1} \left( \zeta^{-1}(X-\zeta)^\ell-\zeta(X-\zeta^{-1})^\ell\right)\ , $$
\noindent
and

$$          {\cal Q}(X) := (\zeta^{-1}-\zeta)^{-1}\left( (X-\zeta)^\ell-(X-\zeta^{-1})^\ell\right) \ .  $$

Note that ${\cal P}(X)$ has degree $\ell$, ${\cal Q}(X)$ has degree $\ell-1$, and both polynomials have coefficients in $K$. It will be convenient to define the rational function $r(X)={\cal P}(X)/{\cal Q}(X)$. Finally, define
\begin{equation} \label{def:F(X,T)}
        F(X,u) ={\cal P}(X)-u{\cal Q}(X)\in K[X,u]\ . 
\end{equation}
Here we assume  $u$ is transcendental over $K$. This is a higher degree analogue of the Shank's polynomial as becomes clear from the following theorem.

\begin{thm}
The polynomial $F(X,u)$ is irreducible over $K(u)$. Let $x$ be a root in some extension field of $K(u)$. Then, $K(x,u)=K(x)$ is a Galois extension of $K(u)$ with Galois group isomorphic to $\mathbb{Z}/\ell \mathbb{Z} $.
                    
                    The discriminant of $F(X,u)$ is given by
                    \begin{equation} \label{e:disc}
                            \ell^\ell(4-\omega^2)^{(\ell-1)(\ell-2)/2}(u^2-\omega u+1)^{\ell-1} \ . \end{equation}
                    \end{thm}
 
 Note that if $x$ is a root of $F(X,u)=0$, then $u={\cal P}(x)/{\cal Q}(x)=r(x)$. This justifies the equality $K(x,u)=K(x)$.
 The formula for the discriminant is stated in Rikuna's paper, but not proven there. A proof can be found in Komatsu [17], Lemma 2.1.
 
 Finally, we note that the polynomial $P(u) = u^2-\omega u+1=(u-\zeta)(u-\zeta^{-1})$ plays a big role in our considerations. From now on we will assume that $\zeta\notin K$. This implies that $P(u)$ is irreducible over $K$. The formula for the discriminant then shows that the only primes of $K(u)$ which can ramify in $K(x)$ are the zero divisor of $P(u)$ and possibly the prime at infinity. A simple calculation,  
 using the Riemann-Hurwitz formula, shows the prime at infinity does not ramify. Thus, $K(x)/K(u)$ ramifies at exactly one prime, the zero divisor of $P(u)$.

\section{Preliminaries}

The following lemma is well-known, and a proof can be found in \cite{Lang}.

\begin{lemma}\label{l:lang} Let $k$ be a field, $m$ an integer $\geq 2$, and $a \in k$, $a \not= 0$. Assume that for any prime $p$ with $p \mid m$, we have $a \notin k^p$, and if $4 \mid m$, then $a \not\in -4k^4$. Then $x^m - a$ is irreducible in $k[x]$.
\end{lemma}

We will also need the following.

\begin{lemma} \label{l:powers}
Let $A$ be an abelian group, and $a$ an element of $A$. Suppose that $a$ is an $n_1$-power and an $n_2$-power with $(n_1,n_2)=1$. Then, $a$ is an $n_1n_2$-power.
\end{lemma}

\noindent
{\it Proof.}  By hypothesis, there exist $b,c\in A$ such that $a=b^{n_1}$ and $a=c^{n_2}$. 
Since $(n_1,n_2)=1$ there exist integers $r$ and $s$ such that $rn_1+sn_2=1$. Then,

$$   a=a^1=(a^r)^{n_1}(a^s)^{n_2}=c^{n_2rn_1}b^{n_1sn_2}=(c^rb^s)^{n_1n_2} \quad. $$ \qed

The main goal of this section is to prove the following.

\begin{lemma} \label{l:gamma}
Let $\ell$ be an odd prime,  $m>1$  an integer not divisible by $\ell$, and $\zeta$ a primitive $l$-th root of unity.  For all sufficiently large prime powers $q$ satisfying 

(i) $q\equiv -1 \pmod{\ell}$, and

(ii) $q\equiv 1\pmod{m_0}$ where $m_0$ is the square-free part of $m$,

 there is a $\gamma\in \F_q^\times$ such that ${X^m-(\gamma+\ell \zeta)}$ is irreducible over $\F_q(\zeta)$.
\end{lemma}

\noindent
{\it Proof}.    We begin by reducing the problem to one which takes place entirely in the field $\F_q$.

Since $q\equiv -1\pmod{\ell}$ it follows that the quadratic extension of $\F_q$ has the form $\F_q(\zeta)$, where $\zeta$ is a primitive $\ell$-th root of unity.  Note that since $\F_q(\zeta) = \F_{q^2}$, then -1 must be a square in $\F_q(\zeta)$.  As a result, to prove that ${X^m-(\gamma+\ell \zeta)}$ is irreducible over $\F_q(\zeta)$, it is enough by Lemma~\ref{l:lang} to show that $\gamma+\ell \zeta$ is not a $p$-th power for all primes $p$ dividing $m$:
 if $\ell\zeta + \gamma = -4w^4$ for some $w \in \F_q(\zeta)$, then $\ell\zeta + 
\gamma = (2w^2\alpha)^2$ is a square in $\F_q(\zeta)$, a contradiction.

So let $p$ be a prime dividing $m$ and suppose that $\gamma + \ell \zeta$ is a $p$-th power in $\F_q(\zeta)$. Taking norms from $\F_q(\zeta)$ to $\F_q$, we find that $\gamma^2+\ell (\zeta+\zeta^{-1})\gamma +\ell^2$ is a $p$-th power in $\F_q$.  Completing the square, we find $c$ and $d$ in $\F_q$ such that $$\gamma^2+(\zeta+\zeta^{-1})\ell\gamma + \ell^2= (\gamma-c)^2+d\ . $$
A short computation shows that $d\ne 0$.  It follows that if we can find a $\gamma \in \F_q$ such that $(\gamma-c)^2+d$ is not a $p$-th power in $\F_q$ for every prime $p|m$, then $X^m-(\gamma+\ell \zeta)$ is irreducible over $\F_q(\zeta)$ as required. We will show that for $q$ large enough there exists $\lambda\in \F_q$ such that $\lambda^2+d$ is not a $p$-th power for every prime $p$ dividing $m$. Then, $\gamma=\lambda+c$ will be the element we are looking for.\par

For each $k$ dividing $q-1$, consider the curve $C_k: y^2+d=x^k$. This curve is absolutely irreducible and non-singular except for the unique point at infinity when $k>3$. Its genus is $(k-1)/2$ when $k$ is odd, and ${k\over 2}-1$ when $k$ is even. Let $N_k$ be the number of points $(\alpha,\beta)\in \F_q^{(2)}$ such that $\beta^2+d=\alpha^k$, i.e. the number of rational points on $C_k$. Using either the Riemann hypothesis for curves, or a more elementary argument using Jacobi sums (see~\cite{IR}, Chapter 8), one can show that $|N_k-q|\le (k-1)\sqrt{q}$. We will need this estimate, especially when $k$ is square-free dividing $m$. Our hypothesis ensures that in this case, $k$ divides $q-1$.\par

Let $R_k$ denote the set of $k$-th powers in $\F_q$ (including zero), and let $$S_k = \{\eta\in R_2\ |\ \eta +d\in R_k\}.$$ It is easy to see that $R_2$ has ${q+1\over 2}$ elements. What can be said about the size of $S_k$? Well, if $(\alpha,\beta)$ is a rational point on $C_k$, i.e. an element of $C_k(\F_q)$,  then $\beta^2\in S_k$. So, there is a map $(\alpha,\beta)\to \beta^2$ from $C_k(\F_q)$ to $S_k$. From the definition of $S_k$, it is clear that this map is onto. Since $\pm 1\in \F_q$ and the $k$-th roots of unity are in $\F_q$, the map is $2k$ to $1$ at all but at most two elements of $S_k$, namely $0$ and $-d$ ($0$ if $d$ is a $k$-th power, and $-d$ if $-d$ is a square). In all cases, one can show that $|\#(S_k)-N_k/2k|< 2$. It follows that the number of elements in $S_k$ is approximately $q/2k$.
\par
If $S$ is a subset of $R_2$, let $S'$ denote its complement in $R_2$. Consider the set

$$                             T=\bigcap_{p|m} {S_p}' \ . $$

\noindent
The intersection is over all primes dividing $m$. If $\tau\in T$, then $\tau+d$ is not a $p$-th power for any prime $p$ dividing $m$. Thus, if $\tau=\lambda^2$ then $\gamma=\lambda +c$ is the element we are looking for. We will show that $T$ is non-empty for $q$ large enough. In fact, we will show a lot more, namely

$$                        \#T = {q\over 2}\prod_{p|m}(1-{1\over p}) + O(\sqrt{q})\ . $$

To this end, let's enumerate the primes dividing $m$, i.e. $p_1,p_2,\dots,p_t$. Then,

$$                     T'=\bigcup_{i=1}^t S_{p_i} \ , $$
\noindent
and therefore,

$$              \#(T') = \sum_i \#(S_{p_i}) -\sum_{i<j} \#(S_{p_i} \cap S_{p_j})+\sum_{i<j<k} \#(S_{p_i}\cap S_{p_j} \cap S_{p_k})- \dots $$
\noindent
by the inclusion/exclusion principle. \par

The intersections simplify considerably. Namely, it can be shown via Lemma~\ref{l:powers} that

$$                    S_{p_{i_1}}\cap S_{p_{i_2}}\cap \dots \cap S_{p_{i_r}} = S_{p_{i_1}p_{i_2} \dots p_{i_r}} \ . $$
\noindent
Since, by hypothesis, the square-free part of $m$ divides $q-1$ we can use our previous estimates, $|\#(S_k)-N_k/2k|<2$ and $|N_k-q|\le k\sqrt{q}$. From this we see

$$                                   \#(S_k)= {q\over 2k} + O(\sqrt {q}) \ , $$

\noindent
for all square-free $k$ dividing $m$. Using this in the above expression for $\#(T')$ yields

$$        2 \#(T')/q = \sum_i {1\over p_i} -\sum _{i<j} {1\over p_ip_j} + \sum_{i<j<k} {1\over p_ip_jp_k} -\dots +O(q^{-{1\over 2}}) \ . $$
\noindent
which is equivalent to (using $\#(R_2)={q+1\over 2}$)

$$          \#(T) ={q\over 2}\prod_i(1-{1\over p_i}) + O(\sqrt{q}) \ . $$

By paying more attention to detail it is fairly easy to give an explicit lower bound for $\#(T)$ in terms of $q$ and thus determine how large $q$ has to be in order to ensure the $T$ is non-empty. See the appendix for details.  \qed

\section{Ichimura's Lemma and Class Number Indivisibility}

In \cite{Ichimura}, Ichimura states a version of the following lemmas, though his proof seems incomplete.  Here we give a rigorous proof, using the same ideas which Iwasawa used in his original result for number fields.

\begin{prop} (Ichimura's Lemma)  \label{l:ichimura}
Let $K/k$ be a finite, geometric, $\ell$-extension which is ramified at exactly one prime $\mathfrak p$ of $k$. Suppose that only one prime $\mathfrak P$ of $K$ lies above $\mathfrak p$, and $\ell \nmid \deg \mathfrak p$. Then, $\ell\ |\  h_K$ implies $\ell\ |\  h_k$. 
 \end{prop}

 First, we fix some notation. Let $k$ be a function field in one variable with finite field of constants $\F_q$. Let $\mathfrak p$ be a prime of $k$ and $A$ the subring of $k$ consisting of elements whose only poles are at $\mathfrak p$. It is well known that $A$ is a Dedekind domain and that its group of units is precisely $\F_q^\times$.\par
 The proof of the following lemma is given in~\cite{Ro1}.
 
\begin{lemma} \label{l:exact_seq}
Let $J_k$ be the group of divisor classes of degree $0$ of k, $Cl_A$ the ideal class group of $A$, and $d=\deg \mathfrak p$. Then, the following sequence is exact.
 $$  (0)\rightarrow J_k\rightarrow Cl_A \rightarrow \Z/d\Z\rightarrow (0)   .  $$
\end{lemma} 

\begin{cor}
Let $h_A=\# Cl_A$, the class number of $A$, and $h_k=\# J_k$, the class number of k. Then
 $$                     h_A=h_k d .   $$
\end{cor}

  A proof of the following can be found  in~\cite{Ro1}.
 
\begin{prop}
Let $k_A$ be the maximal, abelian, unramified extension of $k$ in which $\mathfrak p$ splits completely.  Then $k_A$ is a finite abelian extension of $k$ and
     $$               {\rm Gal}(k_A/k) \cong Cl_A  .  $$
\end{prop}

\smallskip\noindent
{\it Proof of Ichimura's Lemma.}   \quad Let $B$ be the integral closure of $A$ in $K$. Applying Lemma~\ref{l:exact_seq} and its corollary to the pair $B, \mathfrak P$, we see that $\ell\ |\  h_K$ implies $\ell \ |\ h_B$.  Let $E$ be the maximal abelian, unramified, $\ell$-extension of $K$ in which $\mathfrak P$ splits completely. Since $E\subset K_B$, and $\ell \ |\ h_B=[K_B:K]$, we see that E properly contains $K$. \par
It is easily seen that $E/k$ is a Galois $\ell$-extension. Let $G$ denote its Galois group. For a prime $\cal P$ of $E$ lying over $\mathfrak P$, let $D(\cal P/\mathfrak p)$ be its decomposition group over $k$. Note that
$$  |D({\cal P}/{\mathfrak p})| = e({\cal P}/{\mathfrak p}) f({\cal P}/{\mathfrak p})= e({\mathfrak P}/{\mathfrak p}) f({\mathfrak P}/{\mathfrak p})=[K:k]\ . $$
The last inequality is because of the assumption that $\mathfrak P$ is the only prime of $K$ lying over $\mathfrak p$. We conclude that $D(\cal P/\mathfrak p)$ is a proper subgroup of $G$. Since $G$ is an $\ell$-group, it follows from a well known result about $\ell$-groups that $D(\cal P/\mathfrak p)$ is contained in a normal subgroup $N\subset G$ of index $\ell$. Any other prime $\cal P'$ of $E$ over $\mathfrak P$ has a decomposition group over $k$ which is conjugate to $D(\cal P/\mathfrak p)$ and is thus also contained in $N$. It follows that the fixed field $L$ of $N$ is a cyclic, unramified extension of $k$ in which $\mathfrak p$ splits completely. It follows that $L\subset k_A$. Thus, $l\ |\ h_A=h_k d$ by the corollary to Lemma~\ref{l:exact_seq}. Since we are assuming that $\ell$ does not divide $d$, we must have $\ell \ |\ h_k$, as asserted. \qed

We now use Ichimura's lemma to prove the following result.

\begin{thm}  \label{t:ichimura_app}
Let $k/\F_q$ be a function field on one variable of a finite constant field $\F_q$ with $q$ elements. Let $\ell$ be a fixed rational prime, and suppose that $q-1$ is not divisible by $\ell$. Suppose further that the class number $h_k$ of $k$ is not divisible by $\ell$. Then, for every postive integer $t$ there are infinitely many non-isomorphic geometric  extensions $K$ of $k$ such that $[K:k]=\ell^t$ and for which $h_K$ is not divisible by $\ell$.
\end{thm}

A variant of this theorem also holds in the number field case, but we will not prove it here.

\medskip
Before proceeding to the proof, we need to recall some facts about the class field theory of a global function field $k$. We will use the language of valuations rather than primes. As is well known, these are completely equivalent concepts. Let ${\it M}_k$ denote the set of normalized valuations $v$ of $k$. For each $v\in{\it M}_k$, let $k_v$ be the completion of $k$ at $v$, $O_v=\{a\in k_v\ |\ {\rm ord}_v(a)\ge 0\}$, $P_v=\{a\in k_v\ |\ {\rm  ord}_v(a)>0\}$, and $U_v=\{a\in k_v\ |\ {\rm ord}_v(a)=0\}$. Note that $O_v$ is a complete, discrete valuation ring, and $U_v$ is its group of units. The norm of $v$, $Nv$, is, by definition, the number of elements in the residue class field of $O_v$, $\kappa_v:=O_v/P_v$. The ideles of $k$, $I_k$, are the subgroup of the direct product $\prod_v k_v^*$ consisting of elements $(a_v)$ where all but finitely many $a_v\in U_v$. The group $k^*$ injects into $I_v$ on the diagonal. For a fixed $w\in {\it M}_k$, let $U_w^{(1)} = \{a\in U_w\ |\ a\equiv 1\pmod {P_w}\}$.  Finally, define
$$        {\cal U}_k=\prod_v U_v \quad {\rm and}\quad {\cal U}_k(w)=\prod_{v\ne w} U_v \times U_w^{(1)}\ . $$
We have an exact sequence
\begin{equation} \label{E:exact_2}
 (0)\rightarrow k^*{\cal U}_k/k^*{\cal U}_k(w)\rightarrow I_k/ k^* {\cal U}_k(w)\rightarrow I_k/k^*{\cal U}_k \rightarrow (0) \end{equation}

By the main theorem of global class field theory (see Artin-Tate or J. Neukirch), each term in this exact sequence is isomorphic to the Galois group of an arithmetically defined field extension. 
The fourth term is isomorphic to a dense subgroup of the Galois group of $k^{nr}/k$, the maximal, abelian, unramified extension of $k$. Let $\bar k$ be the maximal constant field extension of $k$. We know that ${\rm Gal}(\bar k/k)\cong \hat{\Z}$, and ${\rm Gal}(k^{nr}/\bar k)\cong J_k$. This parallels the exact sequence
$$                       (0)\rightarrow J_k\cong I_k^0/k^*{\cal U}_k\rightarrow I_k/k^*{\cal U}_k\rightarrow \Z\rightarrow (0)\ . $$
The arrow from the third to the fourth term is induced by the degree map on ideles, $\deg: I_k\rightarrow \Z$, given by $(a_v)\rightarrow \sum_v {\rm ord}_v(a_v)\deg v$. Recall that the degree of $v$ is the dimension of $\kappa_v := O_v/P_v$ over $\F_q$.\par
Now consider the third term of the exact sequence~(\ref{E:exact_2}). It is isomorphic to a dense subgroup of the Galois group of $k(w)/k$, the maximal, abelian extension of $k$ which is unramified at all $v\ne w$ and is at most tamely ramified at $w$. Clearly, $k^{nr}\subset k(w)$ and we have
$$  {\rm Gal}(k(w)/k^{nr})\cong k^*{\cal U}_k/k^*{\cal U}_k(w)\cong {\cal U}_k/{\cal U}_k\cap k^*{\cal U}_k(w) =  $$
$${\cal U}_k/\F_q^\times{\cal U}_k(w)\cong U_w/\F_q^\times U_w^{(1)}\cong \kappa_w^*/\F_q^\times\ .$$
\noindent
The last isomorphism comes from the natural reduction map of $U_w$ to $\kappa_w^*$ which is onto with kernel $U_w^{(1)}$.

We have proved the following important lemma.

\begin{lemma} \label{l:Gal}
The Galois group of $k(w)/k^{nr}$ is cyclic of order ${Nw -1\over q-1}$.
\end{lemma}

\noindent
{\it Proof of Theorem \ref{t:ichimura_app}.}  It is enough to prove the the theorem in the case $t=1$. If $K_1/k$ is a geometric, cyclic extension of degree $\ell$ with the property that $\ell$ does not divide $h_{K_1}$, simply replace $k$ by $K_1$ and use the theorem again. In finitely many steps, a field of degree $\ell^t$ will be constructed whose class number is indivisible by $\ell$. It will follow from the construction to be given below that this process will produce infinitely many non-isomorphic fields with the required properties.\par
The field $k(w)$, which we defined in the discussion preceding Lemma~\ref{l:Gal}, is infinite dimensional over $k$. To deal with this we choose an auxiliary valuation $v_o\ne w$ and define $k(w,v_o)$ to be the maximal abelian extension of $k$ which is at most tamely ramified at $w$, in which $v_o$ splits completely, and is unramified for all $v\ne w$. Let $k_o$ be the maximal constant field extension of $k$ in $k(w,v_o)$. Exactly as in the proof of Theorem 1.3 of~\cite{Ro1}, one can show that 
\begin{equation}     [k_o:k]=\deg v_o:=d_o \quad {\rm and}\quad {\rm Gal}(k(w,v_o)/k_o)\cong {\rm Gal}(k(w)/\bar k) \end{equation}
We'll come back to this in a moment.\par
Let $e$ be the order of $q$ mod $\ell$. Since we have assumed $\ell$ does not divide $q-1$ we see that $1<e < \ell $. Let $n$ be any positive integer indivisible by $\ell$ and $w$ a valuation of degree $ne$. If $ne$ is large enough, such a valuation must exist (see Theorem 5.12 in [Ro2]). Note that $\deg w = ne$ is not divisible by $\ell$. \par
Recall that $Nw=q^{\deg w}=q^{ne}$. We claim that $Nw-1/ (q-1)$ is divisible by $\ell$. This follows easily from

\begin{equation}      {Nw-1\over q-1}={q^{ne}-1\over q-1}= {q^{ne}-1\over q^e-1} {q^e-1\over q-1}\ . \end{equation}
\noindent
Both factors are in $\Z$ and the last factor is divisible by $\ell$ since $\ell$ does not divide $q-1$. \par
We now return to equation $(4)$. We have a tower $\bar k\subset k^{nr}\subset k(w)$, and thus a corresponding tower between $k_o$ and $k(w,v_o)$. What is the middle of this tower? Let $A\subset k$ be the ring of elements whose valuations are non-negative away from $v_o$. Then, $k_A$ is the maximal abelian, unramified extension of $k$ in which $v_o$ splits completely. A moments reflection reveals that this is the middle term, i.e. we have $k_o\subset k_A\subset k(w,v_o)$. \par
Having chosen $w$ to have degree $ne$ we now impose on $v_o$ the requirement that its degree $d_o$ be prime to $\ell$. Again, using Theorem 5.12 in [Ro2], we see that such $v_o$ exist in abundance. We  claim that $k(w,v_o)$ contains an intermediate extension $K$ which is geometric, cyclic of degree $\ell$ over $k$, and ramified at $w$ and nowhere else. Since $\deg w=ne$ is indivisible by $\ell$, we can apply Ichimura's result, Proposition $2$, to conclude that the class number of $K$ is not divisible by $\ell$. This will conclude the proof.\par
To find such a $K$, recall (equation $(4)$) ${\rm Gal}(k(w,v_o)/k_o)\cong {\rm Gal}(k(w)/\bar k)$ which has ${\rm Gal}(k(w)/k^{nr})$ as a quotient. By Lemma 3, the order of this Galois group is $Nw-1/(q-1)$. Since $\deg w =ne$, equation $(5)$ shows that $\ell$ divides this number. Consequently, $\ell$ must divide the order of the Galois group of $k(w,v_o)/k$. Since this group is abelian, it must have a subgroup of index $\ell$. Let $K$ be the fixed field of such a subgroup. We claim that $K$ has all the necessary properties.\par
First, we clearly have ${\rm Gal}(K/k)$ is cyclic of order $\ell$. Next, we show that $K/k$ is not a constant field extension. If it were, then $K\subset k_o$. However, $[k_o:k]=d_o$, which we chose prime to $\ell$. Thus, $K/k$ is a geometric extension. Finally, we claim that $K/k$ is ramified at $w$ and nowhere else. Since $K\subset k(w,v_o)$ it is unramified at every valuation $v\ne w$. If it were also unramified at $w$ it would follow that $K\subset k_A$. However, by Proposition $1$ and the Corolllary to Lemma $1$, $[k_A:k] = h_k d_o$ which is prime to $\ell$. Thus, $K/k$ cannot be unramified, so it must be ramified and totally ramified at $w$. We have shown that $K$ satisfies all the necessary properties, so the proof of the theorem is complete. \qed


\section{Proofs of Main Results}

We are now ready to prove Theorem~\ref{t:main_result}.  Let $m$ be any positive integer $m > 1$ and $\ell$ an odd prime. Write $m = \ell^tm_1$ for integers $t$ and $m_1$ with $\ell \nmid m_1$.  Let $m_0$ be the square-free part of $m_1$, and fix a prime power $q$, sufficiently large, with $q \equiv 1$ (mod $m_0$) and $q \equiv -1$ (mod $\ell$).  First, we prove the theorem for the case when $\ell \nmid m$.

Define rational functions $X_j(T)$ recursively as follows,  $X_0(T)=T$ and

\begin{equation} \label{X_j def}
 X_j={{\cal P}(X_{j-1})\over {\cal Q}(X_{j-1})}=r(X_{j-1}) \quad {\rm for}\quad j\ge 1\ . 
 \end{equation}

The relevant notations were introduced at the end of the introduction. Note that $X_j=r^{(j)}(T)$, where the superscript $(j)$ means to compose $r$ with itself $j$ times.\par
Recalling the Rikuna  polynomial $F(X,u) = {\cal P}(X)-u{\cal Q}(X)$ we see that $F(X_{j-1},X_j)=0$. It follows from Theorem 1.4, and the following remarks, that  $\F_q(X_{j-1})/\F_q(X_j)$ is a cyclic extension of degree $\ell$, ramified only at the zero divisor of $X_j^2-\omega X_j+1$.\par
Now, fix a positive integer $n\ge 1$, and for $1\le i \le n$ define 

$$                N_i= \F_q(X_{n-i}) \quad {\rm  and}\quad  M_i=N_i(\root m\of {\ell X_n+\gamma})\ . $$
\noindent
Here $\gamma\in  \F_q$ is chosen so that $X^m-(\ell \zeta + \gamma)$ is irreducible over $\F_q(\zeta)$. See Lemma 2.3.\par
Note that $N_n=\F_q(T)$ and $M_n= \F_q(T)(\root m\of {\ell X_n+\gamma})$ We will show that $M_n$ is an extension of $\F_q(T)$ of degree $m$ and that its class number is not divisible by $\ell$. Further, the genus of $M_n$ is a monotone increasing function of $n$. Thus, all the fields $M_n$ are pairwise non-isomorphic. This will prove our theorem in the case $m$ is not divisible by $\ell$.

We will see that for all $i$ such that $1 \leq i \leq n - 1$, $[N_{i+1} : N_i] = \ell$, $[M_{i+1} : M_i] = \ell$, and for all $i$, $[M_i : N_i] = m$.  The field diagram is shown below.

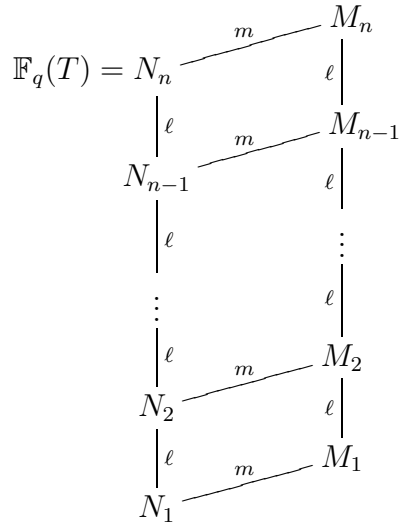
\[\xymatrix@=20pt@R=3pt{
	&**[r]M_n  \ar@{-}[dd]_\ell\\
	**[l]\F_q(T) = N_n \ar@{-}[ur]^m \ar@{-}[dd]^\ell\\
	&**[r]M_{n-1} \ar@{-}[dd]_\ell\\
	N_{n-1} \ar@{-}[ur]^m \ar@{-}[dd]^\ell\\
	&{\vdots} \ar@{-}[dd]_\ell\\
	{\vdots} \ar@{-}[dd]^\ell\\
	&M_2 \ar@{-}[dd]_\ell\\
	N_2 \ar@{-}[ur]^m \ar@{-}[dd]^\ell\\
	&M_1\\
	N_1 \ar@{-}[ur]^m
	}\]

Let
\[P_i = X_{n-i}^2 - \omega X_{n-i} + 1,\]
and let $(P_i)$ denote the divisor of $N_i$ corresponding to the zeros of $P_i$.  Recall that $q \equiv -1$ (mod $\ell$) which implies that $X^2 - \omega X + 1$ is irreducible over $\F_q$. Therefore, $P_i$ is irreducible in $\F_\q[X_{n-i}]$, and hence $(P_i)$ is a prime divisor.

The idea of the proof of the main result is as follows.  We will show that $\ell \nmid h_{M_1}$, and use Proposition~\ref{l:ichimura}  to conclude that $\ell \nmid h_{M_n}$.  The next few lemmas
show that Proposition~\ref{l:ichimura} applies.  Finally, we show that the $M_n$'s are distinct, so there are infinitely many degree $m$ extensions of $\F_\q$ with class number indivisible by $\ell$.

\begin{lemma}\label{N_i+1 cyclic over N_i} For each $i$, $N_{i+1}$ is a $\Z/\ell \Z$-extension of $N_i$, totally ramified at $(P_i)$, and unramified outside $(P_i)$.
\end{lemma}

\begin{proof} By the remarks on the previous page, we see that $N_{i+1}$ is a $\Z/\ell \Z$-extension of $N_i$.  By Eq.(\ref{e:disc}), the discriminant is 
$$
\ell^\ell (4-\omega^2)^{(\ell-1)(\ell-2)/2}(X_{n-i}^2-\omega X_{n-i} + 1)^{\ell-1} = \ell^\ell (4-\omega^2)^{(\ell-1)(\ell-2)/2}P_i^{\ell-1},
$$ where $\ell^\ell (4-\omega^2)^{(\ell-1)(\ell-2)/2} \in \F_\q^\times$. It is easy to see $\ell^\ell (4-\omega^2)^{(\ell-1)(\ell-2)/2} \neq 0$ since $\chr \F_\q \neq \ell$ and if $4 - \w^2 = 0$, then $\w = \pm 2$.  This implies that $\zeta + \zeta^{-1} = \pm 2
$, so $\zeta = \pm 1,$ a contradiction since $\ell \geq 3$.

Since any finite ramified prime would divide the discriminant, it follows that the only possible ramification is at $P_i$ and at the prime at infinity.  Note that the infinite prime has degree 1, so if $(P_i)$ were unramified, then Riemann-Hurwitz implies that
\[2g_{N_{i+1}} - 2 = \ell(2g_{N_i} - 2) + e_{\infty} - 1.\]
Since $N_i$ and $N_{i+1}$ are rational function fields, they both have genus $0$.  It follows that $e_\infty = 2\ell - 1$, which is impossible since the ramification index is at most the degree of the extension, which is $\ell$ in this case.  So $(P_i)$ must be ramified in $N_{i+1}$, and the ramification index is $\ell$ since the extension is Galois of prime degree $\ell$. It follows that the infinite prime is unramified, because
\[-2 = -2\ell + (\ell - 1)\deg(P_i) + e_{\infty} - 1 = -2\ell + 2\ell - 2 + e_{\infty} - 1 = e_{\infty} - 3.\]
So $e_{\infty}=1$, as claimed.
\end{proof}

\begin{lemma}\label{P_1 inert in M_1} The extension $M_i/N_i$ has degree $m$, and the prime $(P_i)$ of $N_i$ is inert in the extension $M_i$.
\end{lemma}

\begin{proof} Since $M_i = N_i(\sqrt[m]{\ell X_n + \gamma})$, it suffices to show that the minimal polynomial for $\sqrt[m]{\ell X_n + \gamma}$ over $N_i$ is irreducible mod $P_i$. We will show that $X^m - (\ell X_n + \gamma)$ is irreducible mod $P_i$, which implies that $X^m - (\ell X_n + \gamma)$ is irreducible over $N_i$ and thus must be the minimal polynomial for $\sqrt[m]{\ell X_n + \gamma}$ over $N_i$. 


	

Let $\lambda$ be the unique $\F_q$ - homomorphism from $\F_q[X_{n-i}]$ to $\F_q(\zeta)$ which takes $X_{n-i}$ to $\zeta$. It is clear that $\lambda$ is onto and has as kernel the principal ideal generated by $P_i$. $\lambda$ extends in the usual way to a homomorphism from the localization $R_i$ of $\F_q[X_{n-i}]$ at the prime ideal $(P_i)$. \par
From definition, we know that $r^{i}(X_{n-i})=X_n$. One easily checks that $r(\zeta)=\zeta$. Using these two facts and $\lambda(X_{n-i})=\zeta$, one deduces that
$\lambda(X_n)=\zeta$. The homomorphism $\lambda$ extends in the obvious way to a homomorphism
from $R_i[X]$ to $\F_q(\zeta)[X]$. This homomorphism takes $X^m-(\ell X_n+\gamma)$ to $X^m-(\ell \zeta +\gamma)$. Since the latter polynomial is irreducible by our choice of $\gamma$, the former one must be irreducible as well. This completes the proof. \end{proof}

\begin{lemma}\label{v(X) separable} The polynomial $\mathcal{Q}(X) \in \F_\q(X)$ is separable.
\end{lemma}

\begin{proof} It suffices to show that $\mathcal{Q}(X)$ and $\mathcal{Q}'(X)$ have no common roots, where $\mathcal{Q}'(X)$ is the formal derivative of $\mathcal{Q}(X)$. The derivative of $\mathcal{Q}(X)$ is given as follows:
\[\mathcal{Q}'(X) = \frac{\ell((X - \zeta)^{\ell-1} - (X - \zeta^{-1})^{\ell-1})}{\zeta^{-1} - \zeta}.\]
Let $\alpha \in \overline{\F_\q}$ be a root of $\mathcal{Q}(X)$. Then, by definition of $\mathcal{Q}(X)$, we have $(\alpha - \zeta)^\ell = (\alpha - \zeta^{-1})^\ell$. Clearly, we cannot have $\alpha = \zeta$ or $\alpha = \zeta^{-1}$, because $\zeta - \zeta^{-1} \not= 0$. If $\alpha$ were also a root of $\mathcal{Q}'(X)$, then we would have $(\alpha - \zeta)^{\ell-1} = (\alpha - \zeta^{-1})^{\ell-1}$.  So \[ (\alpha - \zeta)^\ell = (\alpha - \zeta^{-1})^\ell = (\alpha - \zeta)^{\ell-1}(\alpha - \zeta^{-1}).\]  Since $\alpha \neq \zeta$, then $\alpha - \zeta = \alpha - \zeta^{-1}$, implying that $\zeta = \zeta^{-1}$, a contradiction.
\end{proof}

\begin{lemma}\label{class no. M_1 indiv ell} The class number of $M_1$ is not divisible by $\ell$.
\end{lemma}

\begin{proof} Recall that $M_1 = \F_\q(X_{n-1})(\sqrt[m]{\ell X_n + \gamma})$. First, we claim that the genus of $M_1$ is $(\ell - 1)(m - 1)$. For ease of notation, let $Z = \sqrt[m]{\ell X_n + \gamma}$, so $M_1 = \F_\q(X_{n-1})(Z)$. Notice that $M_1\overline{\F_\q}$ is a degree $m$ extension of $\overline{\F_\q}(X_{n-1})$ with minimal polynomial
\begin{align}\label{min poly of M_1F_q} \notag
	X^m - (\ell X_n + \gamma) & = X^m - \left(\frac{\ell \mathcal{P}(X_{n-1})}{\mathcal{Q}(X_{n-1})} + \gamma \right) \\ \notag
	& = X^m - \frac{\ell \mathcal{P}(X_{n-1})+\gamma \mathcal{Q}(X_{n-1})}{\mathcal{Q}(X_{n-1})} \\
	& = X^m - \frac{F(X_{n-1},-\gamma/\ell)}{\mathcal{Q}(X_{n-1})/\ell}.
\end{align}
(Notice that the polynomial $X^m - (\ell X_n + \gamma)$ remains irreducible over $\overline{\F_q}$:  if $\alpha$ is a zero, then it has multiplicity one; then, in the local ring at $X_{n-1}-\alpha$ the polynomial in question is Eisenstein and so, irreducible.)
The discriminant of $F(X,-\gamma/\ell)$ is $\ell^{-(\ell-2)}(4-\omega^2)^{(\ell-1)(\ell-2)/2}(\gamma^2 + \ell \omega \gamma + \ell^2)$ by Eq.(\ref{e:disc}). This must be non-zero, or else $P_1(-\gamma/\ell) = (\gamma^2 + \omega\gamma\ell + \ell^2)/\ell^2 = 0$. But $-\gamma/\ell \in \F_\q$, and $P_1$ is irreducible over $\F_\q$, a cont
radiction. So $F(X_{n-1},-\gamma/\ell)$ has non-zero discriminant, and hence no multiple roots. By Lemma \ref{v(X) separable}, $\mathcal{Q}(X)$ has no multiple roots.

Finally, $F(X,-\gamma/\ell)$ and $\mathcal{Q}(X)$ must be relatively prime. Otherwise, for some $\alpha \in \overline{\F_\q}$, we would have
\[\mathcal{Q}(\alpha) = 0 = F(\alpha, -\gamma/\ell) = \mathcal{P}(\alpha) + (\gamma/\ell)\mathcal{Q}(\alpha).\]
 It easily follows from the last equality that $\mathcal{P}(\alpha) = 0$. Thus $X - \alpha$ is a common factor of $\mathcal{P}(X)$ and $\mathcal{Q}(X)$ which contradicts the irreducibility of $F(X,u)$.

Hence, the numerator of the constant term in Eq.(\ref{min poly of M_1F_q}) has $\ell$ distinct roots, each corresponding to a prime that is totally ramified in $M_1\overline{\F_\q}$. Similarly, the denominator of the constant term in Eq.(\ref{min poly of M_1F_q}) has $\ell-1$ distinct roots, each corresponding to a prime that is totally ramified in $M_1\overline{\F_\q}$. Finally, it is clear that the infinite prime is totally ramified in $M_1\overline{\F_\q}$. Since $F(X,-\gamma/\ell)$ and $\mathcal{Q}(X)$ are relatively prime, then these $2\ell$ primes are all distinct. Now $\chr \F_\q \nmid m$, and so each of these primes is tamely ramified in $M_1\overline{\F_\q}$. No other primes can be ramified since no other primes can divide the discriminant of $X^m - (\ell X_n + \gamma)$.
Each of the ramified primes has degree 1, so Riemann-Hurwitz implies that
\begin{align*}
	2g_{M_1\overline{\F_\q}} - 2 & = m(2g_{\overline{\F_\q}(X_{n-1})} - 2) + \sum_{\p}(e(\p) - 1)\deg(\p) \\
	& = -2m + 2\ell (m-1) = 2(\ell - 1)(m - 1) - 2,
\end{align*}
and thus $g_{M_1\overline{\F_\q}} = (\ell - 1)(m-1)$, as claimed.

Next, we claim that $M_1 = \F_\q(Z)(X_{n-1})$ is a $\Z/\ell \Z$-extension of $\F_\q(Z)$. We know that $N_1$ is a $\Z/\ell \Z$-extension of $\F_\q(X_n)$ and $\F_\q(Z)$ is a degree $m$ extension of $\F_\q(X_n)$.  (See figure below.)

\[\xymatrix@=7pt@R=10pt{
		& M_1 = \F_q(Z)N_1 \ar@{-}[dl] \ar@{-}[dr]\\
		N_1 = \F_q(X_{n-1}) \ar@{-}[dr]_\ell && \F_q(Z) \ar@{-}[dl]^m\\
		& \F_q(X_n)
	}\]

Since $(\ell, m)=1$, then $M_1 = \F_\q(Z)N_1$ is a $\Z/\ell \Z$-extension of $\F_\q(Z)$. Thus, the minimal polynomial for $X_{n-1}$ over $\F_\q(Z)$ must be $F(X,X_n) = F(X,(Z^m - \gamma)/\ell)$. The discriminant of this polynomial is, by Eq.(\ref{e:disc}),
\[(4-\omega^2)^{(\ell-1)(\ell-2)/2}\ell^\ell (X_n^2 -  \omega X_n + 1)^{\ell-1} = (4-\omega^2)^{(\ell-1)(\ell-2)/2}\ell^{\ell-2(\ell-1)} ((Z^m - \gamma)^2 - \ell \omega (Z^m - \gamma) + \ell^2)^{\ell-1}.\]
Let $(Q)$ be the divisor corresponding to
\[Q = (Z^m - \gamma)^2 - \ell \omega (Z^m - \gamma) + \ell^2 \in \F_\q(Z).\]
We will show that $M_1$ is ramified only at the single prime $(Q)$ of $\F_\q(Z)$, where $\ell \nmid 2m = \deg(Q)$. This completes the proof, by Lemma \ref{l:ichimura}, since $\ell$ does not divide the class number of the rational function field $\F_\q(Z)$. Notice that $Q$ is irreducible over $\F_\q$; if $\alpha$ is a root of $Q$ in some extension of $\F_\q$, then $(\alpha^m - \gamma)/\ell$ is a root of $X^2 - \omega X + 1$, the minimal polynomial of $\zeta^{\pm1}$ over $\F_\q$. So $(\alpha^m - \gamma)/\ell=\zeta^{\pm 1}$.  Since $X^m - (\ell\zeta^{\pm1} + \gamma)$ is irreducible over $\F_\q(\zeta^{\pm1})$, we have $[\F_\q(\alpha):\F_\q(\zeta^{\pm1})] = m$, and so
\[[\F_\q(\alpha) : \F_\q] = [\F_\q(\alpha):\F_\q(\zeta^{\pm1})][\F_\q(\zeta^{\pm1}):\F_\q] = m \cdot 2 = 2m,\]
which proves that $Q$ must be irreducible over $\F_\q$. Thus the divisor $(Q)$ is indeed prime. Since $(Q)$ is the only prime of $\F_\q(Z)$ that divides the discriminant of the minimal polynomial of $X_{n-1}$ over $\F_\q(Z)$, only $(Q)$ and the prime at infinity could be ramified. Assume $(Q)$ is not ramified. By Riemann-Hurwitz, we get
$$2(\ell-1)(m-1) - 2  = (e_{\infty}-1)- 2\ell,$$
so	$e_{\infty} = 2\ell m -2m +1 > \ell$, a contradiction. So $(Q)$ is ramified (totally ramified since the extension is Galois and has prime degree $\ell$) in $M_1$. To see that $M_1$ is ramified at no other primes of $\F_\q(Z)$, we again use the Riemann-Hurwitz formula:
\begin{align*}
	2(\ell-1)(m-1) - 2 & = \ell(-2) + (\ell-1)\deg(Q) + \sum_\p (e_\p - 1)\deg(\p) \\
	& = (\ell - 1)(2m) - 2(\ell - 1) - 2 + \sum_\p (e_\p - 1)\deg(\p) \\
	& = 2(\ell - 1)(m-1) - 2 + \sum_\p (e_\p - 1)\deg(\p).
\end{align*}
Thus, all other primes must be unramified.
\end{proof}

%

\begin{proof}[Proof of Theorem \ref{t:main_result}]
Assume that $\ell \nmid m$.  Notice that $M_{i+1} = M_iN_{i+1}$, so by Lemma \ref{N_i+1 cyclic over N_i}, $M_{i+1}$ is a $\Z/\ell \Z$-extension of $M_i$. Also by Lemma \ref{N_i+1 cyclic over N_i}, $M_{i+1}$ is totally ramified at the prime in $M_i$ lying over $(P_i)$ and unramified everywhere else.  By Lemma~\ref{l:ichimura}, $\ell \nmid h_{M_i}$ implies that $\ell \nmid h_{M_{i+1}}$.  From Lemma~\ref{class no. M_1 indiv ell}, we see that $\ell \nmid h_{M_1}$.  Therefore, $\ell$ does not divide $h_{M_2}, h_{M_3}, ... , h_{M_n}$.  Hence, $M_n$ has class number indivisible by $\ell$.

To show that there are infinitely many such fields, we prove that each $M_n$ has genus $(\ell^n - 1)(m-1)$, so the fields are pairwise non-isomorphic. It was shown in Lemma \ref{class no. M_1 indiv ell} that the genus of $M_1$ is $(\ell-1)(m-1)$.  Since $M_{i+1}/M_i$ is totally ramified at a single prime in $M_i$, denoted here $\mathfrak{P_i}$, lying over $(P_i)$ in $N_i$.  Since $(P_i)$ is inert in $M_i$, $\mathfrak{P_i}$ has degree $2m$ in $M_i$. Note that $M_n$ has degree $\ell^{n-1}$ over $M_1$, so by Riemann-Hurwitz,
\begin{align*}
	2g_{M_n} - 2 & = \ell^{n-1}(2g_{M_1} - 2) + (\ell^{n-1} - 1)(\deg(\mathfrak{P}_1)) \\
	& = \ell^{n-1}(2\ell m - 2\ell - 2m + 2 - 2) + 2\ell^{n-1}m - 2m \\
	& = \ell^{n-1}(2\ell m - 2\ell) - 2m \\
	& = 2\ell^n(m-1) - 2(m-1) - 2 \\
	& = 2(\ell^n - 1)(m-1) - 2.
\end{align*}
Therefore, it follows that $g_{M_n} = (\ell^n - 1)(m-1)$.

Now we consider the general case.  Write $m = \ell^t m_1$, where $\ell \nmid m$, and let $m_0$ be the square-free part of $m_1$.  Since $\ell \nmid m_1$, the results above show that we have infinitely many extensions $K_1$ of degree $m_1$ over $\F_q(T)$ with $\ell \nmid h_{K_1}$.  Note that the constant field of $K_1$ is $\F_q$:  
 $K_1$ is one of the fields $M_n$.
This field is at the top of a tower of totally ramified extensions. At the bottom, $M_1/N_1$ is totally ramified at $X_{n-1}-\alpha$. Also, we know $M_{i+1}/M_i$ is totally ramified at the prime of $M_i$ above $(P_i)$. At a totally ramified prime, the relative degree must be 1. So, in a tower of totally ramified extensions the constant field at the top must be the same as the constant field at the bottom.

 Since $q \equiv -1$ (mod $\ell$), then Theorem~\ref{t:ichimura_app} implies that there are infinitely many non-isomorphic geometric extensions $K$ of degree $\ell^t$ over $K_1$ with $\ell \nmid h_K$. Thus we have infinitely many extensions $K$ of degree $m$ over $\F_q(T)$ with $\ell \nmid h_K$, as claimed. 
\end{proof}

\section{Corollaries}
  
  We are now in a position to prove Corollaries 1.2 and 1.3 which are stated in the introduction. We will reproduce the statements here for the convenience of the reader.

\noindent
{\bf Corollary 1.2}
{\it If $q$ satisfies the hypotheses of the theorem and, in addition, if $q \equiv 1$ (mod $m$), then there are infinitely many cyclic extensions $K$ of degree $m$ over $\F_q(T)$ with $\ell \nmid h_K$.
}

\noindent
{\it Proof}.\quad In the course of the proof of Theorem 1.1, we have produced the following field extensions, $M_n= k(\root m\of{\ell X_n+\gamma}), $ which have degree $m$ and class number indivisible by $\ell$. If $q\equiv 1\pmod{m}$, then the base field contains a primitive $m$-th root of unity. This implies that $M_n$ is a Kummer, and thus cyclic, extension of $k$ of degree $m$. \qed

\noindent
{\bf Corollary 1.3.}   {\it Suppose $m=\ell m_1$ where $\ell$ does not divide $m_1$.  If  $q\equiv 1\pmod{m_1}$,   $q\equiv -1\pmod{\ell}$, then there are infinitely many cyclic extensions $K$ of $k= \F_q(T)$ of degree $m$ with class number indivisible by $\ell$.}

\noindent
{\it Proof}.\quad By Corollary 1.2, there are infinitely many cyclic extensions $K_1$ of $k$ of degree $m_1$ with class number indivisible by $\ell$.  By the proof of Theorem 3.5, there are infinitely many geometric, cyclic extensions $L/k$ of degree $\ell$ whose class number is indivisible by $\ell$. For each such $L$, the compositum, $LK_1$, is a cyclic extension of $k$ of degree $\ell m_1=m$. We will show that for each $K_1$ we can choose an $L$ so that the compositum has class number indivisible by $\ell$. This will prove the corollary.

In the proof of Theorem 3.5, it is shown that for every finite prime $w$ of $k$ of sufficiently large, even degree, there is a cyclic extension $L/k$ which is ramified at $w$ and nowhere else. If the degree of $w$ is not divisible by $\ell$ (which is easy to achieve), then by Ichimura's Lemma, Proposition 3.1, we find the class number of  $L$ must be indivisible by $\ell$. We impose on $w$ one more condition, namely we require the Artin symbol of $w$, $(w,K_1/k)$, be a cyclic generator of ${\rm Gal}(K_1/k)$. By the Chebotarev density theorem there are primes satisfying this condition of all sufficiently large degree (see \cite{Rosen}, Proposition 9.13B). To apply this result it is required that $K_1/k$ be a geometric extension, but this is automatic in this case. Over $\F_q(T)$ any root of a non-constant rational function generates a geometric extension. We omit the elementary proof. 

To summarize, we choose a finite prime $w$ of $k$ of large even degree which is not divisible by $\ell$ and whose Artin symbol $(w,K_1/k)$ is a cyclic generator of ${\rm Gal}(K_1/k)$. Let $L$ be a cyclic extension of $k$ of degree $\ell$ which is ramified at $w$ and nowhere else.\par
The extension $LK_1/K_1$ is cyclic of degree $\ell$ and is ramified only at primes above $w$ in $K_1$. Since $(w,K_1/k)$ is a generator of the Galois group, it follows that $w$ is inert in $K_1$, i.e. there is only one prime $W$ above $w$ and $f(W/w)=m_1$. It follows that the degree of $W$ is $m_1\deg(w)$, and so is not divisible by $\ell$. $W$ is totally ramified in $LK_1$ since $LK_1/K_1$ is a cyclic extension of degree $\ell$. Moreover, $LK_1/K_1$ is unramified at every other prime of $K_1$. If we knew this was a geometric extension, we could invoke Ichimura's Lemma once again to conclude that the class number of $LK_1$ is indivisible by $\ell$. Thus, it only remains to show that $LK_1/K_1$ is a geometric extension.

Let $\mathbb{E}$ be the maximal cofnstant field extension of $\F_q$ in $LK_1$. Since ${\mathbb{E}}\cap K_1=  \F_q$ it follows that $[\mathbb{E} : \F_q]$ is equal to $\ell$ or $1$. However, $\mathbb{E}$ injects into the residue class field of the prime of $LK_1$ lying above $W$. Since $W$ is totally ramified in $LK_1$ we see that $\mathbb{E}$ injects into the residue class field of $W$ which has degree $m_1$ over $\F_q$. Since $\ell$ does not divide $m_1$, we conclude $\mathbb{E} =\F_q$. The corollary is proved.  \qed

\section{Appendix}

The theorem on indivisibility by a prime $\ell$ of the class number of extensions of $\F_q(T)$ of degree $m$ is dependent on the assumption that $q$ is a sufficiently big prime power satisfying
$q\equiv -1\pmod{\ell}$ and $q\equiv 1\pmod{m_o}$, where $m_o$ is the squarefree part of $m$. This is equivalent to a single congruence $q\equiv -1+2\ell \ell' \pmod {\ell m_o}$, where $\ell'$ is a multiplicative inverse of $\ell$ modulo $m_o$. We look into the the question of how big $q$ has to be in order for the theorem to be valid. If $q$ lies in this arithmetic progression and is big enough to make the main theorem valid, we say that $q$ is admissible.\par
The number of rational points on the curve $y^2=x^k-d$ over $\F_q$ satisfied $|N_k-q|\le (k-1)\sqrt{q}$ if $k$ is odd, and $\le 1+(k-1)\sqrt{q}$ if $k$ is even. See Theorem 5 of Chapter 8 in \cite{IR}. The theorem there is stated over the prime field, but the proof work over any finite field. We will work with the slightly weaker, but uniform, inequality $|N_k-q|< k\sqrt{q}$. Also, for the set $S_k$ we have shown
$|\#(S_k)-N_k/2k|<2$. Let's write $N_k=q+\delta_1(k)k\sqrt{q}$ and $\#(S_k) = N_k/2k +2\delta_2(k)$ where $|\delta_1(k)|$ and $|\delta_2(k)|$ are both less than $1$. Putting these two inequalities together, we find

$$                 \#(S_k)= {q\over 2k} +{\delta_1(k)\over 2}\sqrt{q} +2\delta_2(k) \ .  \eqno (1)$$

In the paper, we show that

$$          \#(T')= - \sum_{1<k|m} \mu(k)\#(S_k)\ . $$

Thus, since $\#(T')+\#(T)= (q+1)/2$, we have

$$             \#(T)={q+1\over 2} +\sum_{1<k|m} \mu(k)\#(S_k) \ . \eqno (2) $$

Using equation (1) and substituting into equation (2), yields

$$          \#(T) = {q\over 2}+{1\over 2} +{q\over 2} \sum_{1<k|m} {\mu(k)\over k}+ \sum_{1<k|m} {\mu(k)\delta_1(k)\over 2}\sqrt{q} +2\sum_{1<k|m}\mu(k)\delta_2(k) \ . \eqno (3)$$ 

Combining the first and third terms, simplifies to the following main term

$$                      {q\over 2}\prod_{p|m}(1-{1\over p}) ={q\over 2} {\phi(m_o)\over m_o} \ . $$

To go further, we need the simple observation that $\sum_{k|m}|\mu(k)| =\sum_{r=0}^t {t\choose r}=2^t$ , where $t$ is the number of primes dividing $m$. Since both $\delta_1(k)$ and $\delta_2(k)$ have absolute value less than 1, the sum of the second, fourth, and fifth terms of equation (3) are bounded above by             
$ 2^{t-1}\sqrt{q} + 2^{t+1}$ .  

Putting all this together, we have

$$            | \#(T) - {q\over 2} {\phi(m_o)\over m_o}| \le 2^{t-1}\sqrt{q} + 2^{t+1}\ . $$    

Thus, to insure that $T$ is not empty, it suffices to insure

$$                     q> {2^t m_o\over \phi(m_o)} \sqrt{q}+ 4 {2^t m_o\over \phi(m_o)} \ .  $$

Set $C= 2^tm_o/\phi(m_o)$. The condition can now be written as

$$                 q> C\sqrt{q} +4C \ .                                \eqno (4) $$

Let $f(x)= x^2-Cx-4C$.  The largest zero, $x_o$, of $f(x)$ is given by $2x_o = C+\sqrt{C^2+16C}$. Thus, $x_o$ is less than $C+4$. Equation $(4)$ is satisfied if $f(\sqrt{q})>0$, and this is certainly the case if
$\sqrt{q}> C+4$ since $f(x)$ is easily seen to be increasing at $x_o$ and beyond. We have proved

\begin{prop}
Let $C=2^tm_o/\phi(m_o)$. A prime power $q$ is admissible if $q> (C+4)^2$.
\end{prop}

It is important to point out, that this condition is sufficient but not necessary. We have made a number of somewhat coarse estimates during the derivation. For example, in the case where $\ell =3$ and $m=m_o=2$ (the case considered by Ichimura), every $q$ such that $q\equiv -1\pmod{3}$ is admissible, whereas the Proposition requires $q> 16$. Nevertheless the estimate is strong enough to give some surprising consequences, taking into account the fact that we are looking at $q$ lying in the arithmetic progression $A(\ell,m_o)$ defined by $ q\equiv -1+ll'\pmod{\ell m_o}$. Every $q$ in this progression, except possibly the smallest positive element, is greater than $\ell m_o$. Thus, if $\ell m_o\ge (C+4)^2$, every possible $q$ in this progression with perhaps one exception is admissible. We investigate two special cases.

\begin{cor}      Let's suppose $m_o=p$, a prime. If $p\ge 13$ then every prime power $q$ in $A(\ell,m_o)$ is admissible with at most one exception.
\end{cor}

\noindent
{\it Proof.} If $p\ge 13$ we claim that $\ell p\ge (C+4)^2$ for any odd prime $\ell$. First, let's write out this condition explicitly.

$$                   \ell p \ge ({2p\over p-1}+4)^2 = 4 ({p^2\over  (p-1)^2}+ {4p\over p-1} +4)\ . $$

\noindent
Dividing both sides by $4p$, yields

$$               {\ell\over 4}\ge {p\over (p-1)^2}+{4\over p-1}+{4\over p}\ .                   $$

\noindent
For $p\ge 13$ the right hand side is less than $.74$, so the inequality is satisfied if $\ell$ is greater than $2.96$. Since $\ell$ is an odd prime, this condition is always satisfied.

\begin{cor} 
Suppose $m_o$ is divisible by two or more primes and that the smallest prime dividing $m_o$ is greater than or equal to $7$. Then every prime power $q$ in $A(\ell, m_o)$ is admissible with at most one exception.
\end{cor}

\noindent
{\it Proof.}   \quad The condition we need is

$$                      \ell m_o \ge  16 ({2^{t-2}m_o\over \phi(m_o)} +1)^2 \ . $$

\noindent
Dividing both sides by $16 m_o$ and simplifying yields

$$                   {\ell\over 16} \ge {2^{2t-4}m_o\over \phi(m_o)^2} + {2^{t-1}\over \phi(m_o)} +{1\over m_o}\ . $$

\noindent
If the right hand side of this inequality were less than or equal to $3/16$ this would hold for all odd primes, and the corollary would follow.\par
 An elementary argument shows if $t\ge 2$ the largest value of the right hand side occurs for $m_0=77=7\cdot 11$. In this case the right hand side is

$$                       {77\over 60^2}+{2\over 60}+{1\over 77} \approx  .0677  \ , $$

\noindent
which is comfortably less than $3/16$.

\pagebreak

Allison M. Pacelli,
Williams College,
Williamstown, MA 01267
apacelli@williams.edu

Michael Rosen,
Brown University,
Providence, RI 02906
mrosen@math.brown.edu

\end{document}